\documentclass{amsart}
\usepackage{amsfonts}

\newtheorem{prop}{\large\bf Proposition}
\newtheorem{thm}{\large\bf Theorem}
\newtheorem{dfn}{\large\bf Definition}
\newtheorem{lemma}{\large\bf Lemma}
\newtheorem{Claim}{\large\bf Claim}

\newcommand \pf {\par\noindent{\em Proof:\/} }

\begin{document}

\title{Automorphism Classes of Elements in Finitely Generated
Abelian Groups.}

\author{Charles F. Rocca Jr.}

\curraddr{Department of Mathematics, Western Connecticut State University,  Danbury, CT 06810}

\email{roccac@wcsu.edu}


\maketitle

\begin{abstract}
We will show that every element of a finitely generated abelian
group is automorphically equivalent what we will define to be a {\em representative
element} in a {\em repeat-free subgroup}, and for finite abelian
groups we can count the number of automorphism classes of elements.
\end{abstract}

\section{Representative Elements}

\noindent
 A finite abelian $p$-group $G_p$ can be written
\[G_p=\mathbb{Z}_{p^{r_1}}^{k_1}\oplus\cdots\oplus\mathbb{Z}_{p^{r_n}}^{k_n},\: r_i<r_{i+1}\ for\ all\ i.\]
However for this paper we will write such a group as follows
\[G_p=[\mathbb{Z}_{p^{r_1}}\oplus \cdots \oplus \mathbb{Z}_{p^{r_n}}]
\oplus[\mathbb{Z}_{p^{r_1}}^{k_1-1}\oplus\cdots\oplus\mathbb{Z}_{p^{r_n}}^{k_n-1}],\:
r_i<r_{i+1}.\] This splits the group into two subgroups, a {\em
repeat-free subgroup},
\[G_p^{rf}=\mathbb{Z}_{p^{r_1}}\oplus \cdots \oplus \mathbb{Z}_{p^{r_n}},\]
which contains one copy of each of the factors $\mathbb{Z}_{p^{r_i}}$
and a {\em remainder subgroup},
\[G_p^{rm}=\mathbb{Z}_{p^{r_i}}^{k_1-1}\oplus\cdots\oplus\mathbb{Z}_{p^{r_n}}^{k_n-1},\]
that contains the remaining factors of the group.

For the remainder of this article let $G$ be a finitely generated
abelian group of the form
\[G=G_{p_1}\oplus\cdots\oplus G_{p_m} \oplus \mathbb{Z}^l,\]
where $l\mbox{ and }m \geq 0$ and the $p_i$ are distinct primes
with $p_i<p_{i+1}$ for all $i$.  The above definition of a repeat free
subgroup can be extended to any finitely generated abelian group
$G$ as follows.


\begin{dfn}
Given a finitely generated abelian group $G$ a {\em repeat-free
subgroup} has one of the following forms
\begin{enumerate}
    \item $G^{rf}=G_{p_1}^{rf}\oplus\cdots\oplus G_{p_m}^{rf} \oplus
    \mathbb{Z}$ if $l$ and $m\geq 1$, or
    \item $G^{rf}=G_{p_1}^{rf}\oplus\cdots\oplus G_{p_m}^{rf}$ if $l=0$ and $m\geq 1$, or
    \item $G^{rf}=\mathbb{Z}$ if $l\geq 1$ and $m=0$.
\end{enumerate}
\end{dfn}

\noindent%
The next definition defines our set of {\em representative
elements} in a repeat-free finite abelian $p$-group.

\begin{dfn}
An element $g=(g_1,\ldots,g_n)$ of a repeat-free finite abelian
\mbox{$p$-group}\/ $G_p^{rf}=\mathbb{Z}_{p^{r_1}}\oplus \cdots
\oplus \mathbb{Z}_{p^{r_n}}$ $r_i<r_{i+1}$ for all $i$ is a {\em
representative element} if
\begin{enumerate}
    \item for all\/ $i$ $g_i$ is either $0$ or a power of $p$,
    \item for all\/ $i,j$ if\/ $i<j$, then $g_i<g_j$, and
    \item for all\/ $i,j$ if\/ $i<j$, then the order of $g_i$ in
    $\mathbb{Z}_{p^{r_i}}$ is less than the order of $g_j$ in
    $\mathbb{Z}_{p^{r_j}}$.
\end{enumerate}
\end{dfn}

\noindent%
We will show that in a  finite abelian $p$-group every
element is automorphically equivalent to a representative element of a
repeat-free subgroup.  Therefore we let a {\em representative
element} of a finite abelian $p$-group be an element of a
repeat-free subgroup which is a representative element in that
subgroup. Extending this definition to all finitely
generated abelian groups we get.

\begin{dfn}
Let $G$ ba a finitely generated abelian group, then we say
$g=(g_1,g_2,\ldots,g_m,Z)\in G$ is a {\em representative element}
if
\begin{enumerate}
    \item for all $j$, $g_j$ is a representative element of $G_{p_j}$,
    \item the element $Z\in \mathbb{Z}^l$ is of the form $(z,0,\ldots,0)$ and
    \item for all $j$, each term of $g_j$ is relatively prime to $z$.
\end{enumerate}
\end{dfn}
So, the primary goal of this paper is the following theorem.

\begin{thm}\label{thm:main}Every element of a finitely generated abelian group is
automorphically equivalent to an unique representative element.
\end{thm}

\noindent The proof of this result has two main steps:
\begin{enumerate}
\item Every element is automorphically equivalent to an element in
a repeat-free subgroup.

\item Every element of a repeat-free subgroup is automorphically
equivalent to a representative element.
\end{enumerate}

\noindent The following lemma, in which we think to endomorphisms of finite abelian groups as matrices, is a major tool in proving these
results.

\begin{lemma}{\em \cite{RT}}\label{lemma:autolemma}
An endomorphism $\varphi$ of $G_p=\mathbb{Z}_{p^{r_1}} \oplus
\mathbb{Z}_{p^{r_2}} \oplus \cdots \oplus \mathbb{Z}_{p^{r_m}}$,
$p$ prime and $r_i \leq r_{i+1}$ for all $i$,  is an automorphism
if and only if the reduction of $\varphi$ modulo $p$, written
$\varphi_p$, is an automorphism of $\mathbb{Z}_p^m$.
\end{lemma}

\noindent
The following proposition establishes the first step in
the proof of the theorem.

\begin{prop}\label{prop:prop1}
Every element in a finitely generated abelian group $G$ is
automorphically equivalent to an element in the repeat-free
subgroup.  In particular:
\begin{enumerate}
\item[a)]Each element $g=(g_1,\ldots,g_n)\in\mathbb{Z}_{p^r}^n$,
$n \geq 1$ is automorphically equivalent to an element of the form
$(p^l,0,\cdots,0)$, where the order of $g$ is $p^{r-l}$.
\item[b)]Each element $g=(g_1,\ldots,g_k)$ of $\mathbb{Z}^k$ is
automorphically equivalent to one of the form $(d,0,\cdots,0)$,
where $d=\gcd(g_1,...,g_k)$.
\end{enumerate}
\end{prop}

\pf Let $g\in \mathbb{Z}_{p^r}^n$; each term in $g$ is of the form
$g_i=a_ip^{l_i}$, where $a_i$ is relatively prime to p. Therefore
$g$ is mapped automorphically to the element
$g'=(p^{l_1},\cdots,p^{l_n})$ by multiplying each $g_i$ by the
inverse of $a_i$ modulo ${p^r}$.

The order of $g'$ which is the same as the order of
$g$,  is $p^{r-l}$.  Therefore at least one term of $g'$ is equal to $p^l$,
without loss of generality let $g'_1=p^l$. Also, every term of
$g'$ will be divisible by $p^l$, that is $l_i \geq l$ for all $i$. Therefore
\[
 \left[
 \begin{array}{cccc}
 1 & 0 & \cdots & 0\\
 -p^{l_2-l}&1&\cdots &0\\
 \vdots & \; &\ddots &\vdots\\
 -p^{l_n-l}&0&\cdots &1
 \end{array}
 \right]
 \left[
 \begin{array}{c}
 p^{l}\\
 p^{l_2}\\
 \vdots\\
 p^{l_n}\\
 \end{array}
 \right]=
 \left[
 \begin{array}{c}
 p^{l}\\
 0\\
 \vdots\\
 0\\
 \end{array}
 \right],
\]
where it is clear from the lemma \ref{lemma:autolemma} that this
matrix is an automorphism, thus proving part (a).


\noindent%
The element $g=(g_1,\ldots,g_k)$ of $\mathbb{Z}^k$ can be reduced
using the Euclidian Algorithm by the repeated application of
elementary row operations and this establishes (b).\qed

To show that every element of a repeat-free subgroup is
automorphically equivalent to a representative element we first
restrict ourselves to a finite repeat-free $p$-group and define a
{\em basic reduction} of an element.

\begin{dfn}
Given a repeat-free $p$-group
\[G_p^{rf}=\mathbb{Z}_{p^{r_1}}\oplus\cdots\oplus\mathbb{Z}_{p^{r_n}},\ r_i<r_{i+1}\ for\ all\ i\] and an
element $g=(p^{l_1},\ldots,p^{l_n})$, a {\em Basic Reduction about
position i of $g$\/} is the following:

\par\noindent\hspace*{20pt}
\mbox{for every $j \neq i$ with $p^{l_j} \geq p^{l_i}$ and
$p^{r_j-l_j}$ $\leq$ $p^{r_i-l_i}$ replace $p^{l_j}$ with 0,}

\noindent%
i.e. if the $j^{th}$ term has a greater value but lesser order
than the $i^{th}$ term replace it by 0. \par If the number of
non-zero terms decreases, then the basic reduction is a {\em
non-trivial basic reduction}, otherwise it is a {\em trivial
reduction}. We will say that an element is {\em reduced} if there
are no non-trivial basic reductions.
\end{dfn}

\noindent The matrix for the reduction transformation has elements
\[a_{lk}=\left\{
 \begin{array}{cl}
 1 &if\; l=k\\
 -p^{l_j-l_i} & if\; l=j\ and\  k=i\\
 0 & otherwise
 \end{array} \right.\]
which is similar to the matrix in Proposition \ref{prop:prop1} and is again by by Lemma \ref{lemma:autolemma} is an automorphism.

\begin{lemma}
Let
\[G_p^{rf}=\mathbb{Z}_{p^{r_1}}\oplus\cdots\oplus\mathbb{Z}_{p^{r_n}},\ r_i<r_{i+1}\ for\ all\ i\]
be a repeat-free $p$-group.  Then an element $g \in G_p^{rf}$ is
reduced if and only if $g$ is a representative element.
\end{lemma}

\pf Let $g=(p^{l_1},\ldots,p^{l_n})$ be a reduced element of
$G_p^{rf}$ and let $1 \leq i < j \leq n$.

{\bf Case 1:} If $p^{l_i} \geq p^{l_j}$, then the order of
$p^{l_i}$ in $\mathbb{Z}_{p^{r_i}}$ is strictly less than the
order of $p^{l_j}$ in $\mathbb{Z}_{p^{r_j}}$. Therefore, we can
perform a basic reduction about position $j$, contradicting the
assumption that $g$ is reduced.

{\bf Case 2:} If $p^{l_i} < p^{l_j}$ and the order of $p^{l_i}$ in
$\mathbb{Z}_{p^{r_i}}$ is greater than or equal to the order of
$p^{l_j}$ in $\mathbb{Z}_{p^{r_j}}$, then we can perform a basic
reduction about position $i$ again contradicting the assumption
that $g$ was reduced. Therefore, if $g$ in $G_p^{rf}$ is reduced,
then it is a representative element.

Now suppose that $g$ is a representative element.  Let $1 \leq i <
j \leq n$ so that $p^{l_i} < p^{l_j}$ and the order of $p^{l_i}$
in $\mathbb{Z}_{p^{r_i}}$ is strictly less than the order of
$p^{l_j}$ in $\mathbb{Z}_{p^{r_j}}$, when $p^{l_i}$ and $p^{l_j}$
are non-zero. If we performed a basic reduction about position
$i$, then $p^{l_j}$ would be unchanged since the order of
$p^{l_i}$ is less than that of $p^{l_j}$. Similarly, if we
performed a basic reduction about position $j$, $p^{l_i}$ would
remain unchanged since $p^{l_i} < p^{l_j}$. Since $i$ and $j$ are
arbitrary if $g$ is a representative element, then $g$ is reduced.
\qed

\noindent We are now in a position to prove the main result for
the case of repeat-free $p$-groups.

\pf(Theorem \ref{thm:main} for repeat-free $p$-groups) Let
$g=(g_1,\ldots,g_n)$ be an element in a finite repeat-free
$p$-group $G_p^{rf}$. Applying Proposition \ref{prop:prop1} we know that $g$ is
automorphically equivalent to some $g'$ in which each term is a
power of $p$. Now, if $g'$ is not a reduced element, then we can
perform a non-trivial basic reduction and the the resulting
element, which is automorphically equivalent to $g'$, will have
strictly fewer non-zero terms. Since the rank of $G_p^{rf}$ is
finite the process of making non-trivial reductions will
terminate, and the resulting element will be a reduced element and
therefore a representative element.

Finally we show that representative elements of automorphism
classes are unique.

\begin{Claim}
If $g$ and $h$ are automorphically equivalent representative
elements in $G_p$, then $g=h$.
\end{Claim}
\pf(of claim) Let $g=(g_1,\ldots,g_n)$, $h=(h_1,\ldots,h_n)$ and
let $A$ be the matrix representing the automorphism taking $h$ to
$g$. Suppose $h_i \neq 0$ and let $p^{l_i}$, $0 \leq l_i \leq r_i$
be the order of $h_i$ in $\mathbb{Z}_{p^{r_i}}$ so that we may
write $h_i=p^{r_i-l_i}$.  By Lemma \ref{lemma:autolemma} the
automorphism $A$ must have the form
\[A=
 \left[
 \begin{array}{cccc}
 a_{11}             & a_{12}            & \cdots & a_{1n} \\
 p^{r_2-r_1} a_{21} & a_{22}            & \cdots & a_{2n}\\
 \vdots             &                   &        & \vdots \\
 p^{r_n-r_1} a_{n1} &p^{r_n-r_2} a_{n2} & \cdots & a_{nn}\\
 \end{array}
 \right]
\]

\noindent where $a_{ii}$ is relatively prime to $p$.  Therefore, assuming that $h_i\neq 0$,
\[
 g_i=p^{r_i-l_1} a_{i1}+\cdots+p^{r_i-l_{i}}a_{ii}+\cdots+p^{r_n-l_{n}}a_{in}.
\]

\noindent Since $h$ is a
representative element, for any non-zero term $h_j$ if $j<i$, then
the order of $h_j$ in $\mathbb{Z}_{p^{r_j}}$ is less than the
order of $h_i$ in $\mathbb{Z}_{p^{r_i}}$; i.e. $l_j<l_i$.
Hence, $r_i-l_j>r_i-l_i$ and $p^{r_i-l_i}$ divides
$p^{r_i-l_j}$. If $j>i$, then $h_i<h_j$ and $p^{r_i-l_i}$ divides
$p^{r_j-l_j}$. Thus
\[
 g_i=p^{r_i-l_i}(p^{l_i-l_1} a_{i1}+\cdots+a_{ii}+\cdots+p^{r_n-l_{n}-(r_i-l_i)}a_{ni}).
\]

\noindent Since $a_{ii}$ is relatively prime to $p$
\[
 b_i=(p^{l_i-l_1} a_{i1}+\cdots+a_{ii}+\cdots+p^{r_n-l_{n}-(r_i-l_i)}a_{ni}).
\]

\noindent%
is also relatively prime to $p$. Thus, when $h_i$ is
non-zero, $g_i=b_ih_i$, $b_i$ relatively prime to $p$. However $g$
is a representative element therefore $g_i$ is a power of $p$,
$b_i=1$ and $g_i=h_i$.

An identical argument shows that $h_i=g_i$ when $g_i$ is
non-zero.  Therefore for all $i$, $g_i=h_i$ and so $g=h$. \qed

If $G_p$ is any finite abelian $p$-group, then every element of
$G_p$ is automorphically equivalent to a unique representative
element in $G_p^{rf}$, which is also a representative element for
$G_p$.  Finally, if $G$ is any finitely generated abelian group,
then we have established that any element of $G$ is
automorphically equivalent to element satisfying the first two
conditions of the definition of representative element.  It can be
shown that any element satisfying the first two conditions is
automorphically equivalent to one that satisfies the third
condition and has the same number or fewer nonzero terms. As an
illustration consider the following example.

Let
$G^{rf}=\mathbb{Z}_{p^{r_1}}\oplus\mathbb{Z}_{p^{r_2}}\oplus\mathbb{Z}$
and let $g=[p^{l_1},p^{l_2},kp]$ be any element of $G^{rf}$
satisfying the first two conditions of the definition of a
representative element but not the third.  Then
\[
\left[\begin{array}{ccc}
1 & 0 & p^{r_1-1}-p^{l_1-1} \\
0 & 1 & p^{r_2-1}-p^{l_2-1} \\
0 & 0 & 1
\end{array}\right]
\left[\begin{array}{c}%
p^{l_1} \\
p^{l_2} \\
kp \\
\end{array}\right]=
\left[\begin{array}{c}%
0 \\
0 \\
kp \\
\end{array}\right]
\]
and so $g$ is automorphically equivalent to an element satisfying
the third condition of our definition and has strictly fewer
nonzero terms.


\section{Counting Automorphism Classes}

Having completed the main result we will now show how to count the
number of representative elements in a given finite abelian group.

Let
$G_p^{rf}=\mathbb{Z}_{p^{r_1}}\oplus\mathbb{Z}_{p^{r_2}}\oplus\cdots\oplus\mathbb{Z}_{p^{r_k}}$
where $r_i<r_{i+1}$ for all $i$ be a repeat-free $p$-group and
$g=(p^{l_1},p^{l_2},\ldots,p^{l_k})$ be an automorphism class
representative in $G_p^{rf}$.  The non-zero terms of $g$ are both
increasing and order increasing, thus for $0\leq i<j\leq k$ if
$p^{l_i}$ and $p^{l_j}$ are non-zero terms of $g$ we know that
 \[p^{l_i} < p^{l_j}\]
and
 \[p^{r_i-l_i} < p^{r_j-l_j}.\]
Hence,
 \[l_i < l_j < l_i+(r_j-r_i).\]
We shall refer to the value $r_j-r_i$ as the {\em gap} between the
terms $\mathbb{Z}_{p^{r_i}}$ and $\mathbb{Z}_{p^{r_j}}$ of
$G_p^{rf}$, and in particular we will denote $r_{i+1}-r_i$ by
$n_i$. As an immediate consequence of the above conclusion we know
that if $p^{l_i}$ is non-zero, then either $p^{l_j}$ is zero or we
have at most $r_j-r_i-1$ choices for $l_j$. Further, if $j=i+1$
and $n_i=r_{i+1}-r_i=1$, then one or the other of these two terms
must be zero. Therefore, if the gap between two terms of a group
is one, then in every automorphism class representative at least
one of the terms will be zero.

For motivation let us count the number of automorphism class
representatives in
$G_p^{rf}=\mathbb{Z}_{p^{r_1}}\oplus\mathbb{Z}_{p^{r_2}}\oplus\mathbb{Z}_{p^{r_3}}$
where $r_1<r_2<r_3$. We begin by counting the number of
automorphism classes with a given number of non-zero terms.
\[\begin{array}{|c|l|}\hline
\mbox{\# of Non-Zero Terms} & \mbox{\# of Automorphism Classes}\\
\hline
0 & 1\\
1 & r_1+r_2+r_3\\
2 & r_1(n_{1}-1)+r_1(r_3-r_1-1)+r_2(n_{2}-1)\\
3 & r_1(n_{1}-1)(n_{2}-1)\\\hline
\end{array}\]

\noindent However $r_2=r_1+n_{1}$ and $r_3=r_1+n_{1}+n_{2}$
therefore the above equations can be rewritten as,

\[\begin{array}{|c|l|}\hline
\mbox{\# of Non-Zero Terms} & \mbox{\# of Automorphism Classes}\\
\hline
0 & 1\\
1 & r_1+(r_1+n_{1})+(r_1+n_{1}+n_{2})\\
2 & r_1(n_{1}-1)+r_1(n_{1}+n_{2}-1)+(r_1+n_{1})(n_{2}-1)\\
3 & r_1(n_{1}-1)(n_{2}-1)\\\hline
\end{array}\]

\noindent which gives us,

\[\begin{array}{|c|l|}\hline
\mbox{Non-Zero Terms} & \mbox{\# of Automorphism Classes}\\ \hline
0 & 1\\
1 & 3r_1+2n_{1}+n_{2}\\
2 & r_1(2n_{1}+2n_{2}-3)+n_{1}n_{2}-n_{1}\\
3 & r_1(n_{1}n_{2}-n_{1}-n_{2}+1)\\\hline
\end{array}\]

\noindent and the number of automorphism classes of elements
equals

\[r_1(1+n_{1}+n_{2}+n_{1}n_{2})+1+n_{1}+n_{2}+n_{1}n_{2}=(r_1+1)(n_{1}+1)(n_{2}+1).\]

\begin{prop}
The number of automorphism classes of elements in a repeat-free
finite abelian $p$-group
\[G_p^{rf}=\mathbb{Z}_{p^{r_1}}\oplus\mathbb{Z}_{p^{r_2}}\oplus\cdots\oplus\mathbb{Z}_{p^{r_k}},\: r_i<r_{i+1}\ for\ all\ i\]
is equal to \[(r_1+1)(n_{1}+1)(n_{2}+1)\cdots(n_{k-1}+1),\] where
\[n_{i}=r_{i+1}-r_i.\]
\end{prop}

\pf  We have shown above that this is true for a group with three
terms.  We proceed by induction on the number of terms in
$G_p^{rf}$. By induction the group
\[H_p^{rf}=\mathbb{Z}_{p^{r_1}}\oplus\mathbb{Z}_{p^{r_2}}\oplus\cdots\oplus\mathbb{Z}_{p^{r_{k-1}}}\]
has $(r_1+1)(n_{1}+1)(n_{2}+1)\cdots(n_{k-2}+1)$ automorphism
classes of elements.  This is also the number of automorphism
class representatives in $G_p^{rf}$ in which the $k^{th}$ term is
0.  Therefore we complete the proposition with the following
claim.

\begin{Claim}
If $1<j\leq k$, then in
\[G_p^{rf}=\mathbb{Z}_{p^{r_1}}\oplus\mathbb{Z}_{p^{r_2}}\oplus\cdots\oplus\mathbb{Z}_{p^{r_k}}\]
the number of automorphism class representatives in which the
$j^{th}$ term is the last non-zero term is
\[(r_1+1)(n_{1}+1)(n_{2}+1)\cdots(n_{j-2}+1)n_{j-1}.\]
\end{Claim}

\pf(of claim)\ If $j=2$, then there are $r_2$ automorphism class
representatives in which the second term is the only non-zero term
and $r_1(n_{1}-1)$ in which the first and second terms are the
only non-zero terms.  Therefore the number of automorphism class
representatives in which the second term is the last non-zero term
is

\begin{eqnarray*}
r_1(n_{1}-1)+r_2&=&r_1n_{1}-r_1+r_1+n_{1}\\
                  &=&(r_1+1)n_{1}.
\end{eqnarray*}

Now let $2< j \leq k$, if $2\leq i <j$ and the $i^{th}$ term is
the last non-zero term of an element prior to the $j^{th}$ term,
then there are $(r_j-r_i-1)$ non-zero choices for the value of the
$j^{th}$ term and by induction
\[(r_1+1)(n_{1}+1)(n_{2}+1)\cdots(n_{i-2}+1)n_{i-1}(r_j-r_i-1)\]automorphism class
representatives in which the $j^{th}$ term is the last non-zero
term and the $i^{th}$ term is the second to last non-zero term. So
letting
\[A_i=(r_1+1)(n_{1}+1)(n_{2}+1)\cdots(n_{i-2}+1)n_{i-1}\] and
summing over all $i$ between $2$ and $j$, the total number of
automorphism classes in which the $j^{th}$ term is the last
non-zero term is
\[A_j=r_j+r_1(r_j-r_1-1)+\sum_{i=2}^{j-1}A_i(r_j-r_i-1).\]
However for all $1\leq i < j \leq k$:
\[r_j=r_{j-1}+n_{i-1}\]
and
\[(r_j-r_i-1)=(r_{j-1}-r_i-1+n_{j-1})\] 
so we may rewrite the previous equation as
\[(r_{j-1}+n_{j-1})+(r_1(r_{j-1}-r_1-1)+r_1n_{j-1})+\sum_{i=2}^{j-1}(A_i(r_{j-1}-r_i-1)+A_in_{j-1})\]
And since
\[A_{j-i}=r_{j-1}+r_1(r_{j-1}-r_1-1)+\sum_{i=2}^{j-2}A_i(r_{j-1}-r_i-1)\]
we get

\begin{eqnarray*}
A_{j}&=&A_{j-1}+n_{j-1}+r_1n_{j-1}+A_{j-1}(r_{j-1}-r_{j-1}-1)+\sum_{i=2}^{j-1}A_in_{j-1}\\
&=&A_{j-1}-A_{j-1}+\left(r_1n_{j-1}+n_{j-1}+\sum_{i=2}^{j-1}A_in_{j-1}\right)\\
&=&r_1n_{j-1}+n_{j-1}+\sum_{i=2}^{j-1}A_in_{j-1}\\
&=&n_{j-1}\left((r_1+1)+\sum_{i=2}^{j-1}A_i\right)\\
&=&(r_1+1)(n_{1}+1)(n_{2}+1)\cdots(n_{j-2}+1)n_{j-1},
\end{eqnarray*}
thus proving the claim. In order to finish the proposition we
observe that, by induction, the number of automorphism classes in
which the $k^{th}$ term is zero is
\[(r_1+1)(n_{1}+1)(n_{2}+1)\cdots(n_{k-2}+1)\]
and from the claim the number of automorphism classes in which the
$k^{th}$ term is non-zero is
\[(r_1+1)(n_{1}+1)(n_{2}+1)\cdots(n_{k-2}+1)n_{k-1},\]
and the sum of these two gives the desired result.\qed

Since automorphisms respect the prime decomposition of finite
abelian groups we get this final general result.

\begin{thm}
The number of automorphism classes of elements in a finite
abelian group
\[G=G_{p_1}\oplus\cdots\oplus G_{p_m}\]
where
\[G_{p_j}=\mathbb{Z}_{p_j^{r_1}}^{k_1}\oplus\cdots\oplus\mathbb{Z}_{p_j^{r_n}}^{k_n},\: r_i<r_{i+1}\]
is equal to the product of the number of automorphism classes in
each individual $G_{p_j}^{rf}$.
\end{thm}

\end{document}